\documentclass{amsart}
\usepackage{amssymb,amsmath,latexsym}
\usepackage{amsthm}
\usepackage{fontenc}
\usepackage{amssymb}
\numberwithin{equation}{section}

\newtheorem{theorem}{Theorem}[section]

\begin{document}
\author{Alexander E Patkowski}
\title{On some Integrals associated with the Riesz function}

\maketitle
\begin{abstract}We explore some integrals associated with the Riesz function and establish relations to other functions from number theory that have appeared in the literature. We also comment on properties of these functions.\end{abstract}


\section{Introduction}
A now well-known criterion for the Riemann hypothesis was offered by Riesz [10] (see also [12, pg.382] and [14]), who stated that, a necessary and sufficient condition for the Riemann hypothesis is
\begin{equation}x\sum_{n\ge1}\frac{\mu(n)}{n^2}e^{-x/n^2}=O(x^{\frac{1}{4}+\epsilon}).\end{equation}
Throughout this paper $\mu(n)$ will denote the M$\ddot{o}$bius function [12], and $\dot{R}(x)$ will denote the Riesz series on the left side of (1.1).
In [2, 3] Bartz gave an explicit formula for Merten's function [12, pg.372, Theorem 14.27] $\sum_{n\le x}\mu(n)$ without the assumption of the Riemann hypothesis, and subsequently described the analytic character of two functions defined when there are no multiple zeros (non-trivial)
of the Riemann zeta function. One of these functions is defined by 
$$m(z)=\lim_{n\rightarrow\infty}\sum_{0<\Im\rho<T_n}\frac{e^{\rho z}}{\zeta'(\rho)}.$$
Bartz further offers that [2, Theorem 2], this function can be continued analytically to a meromorphic function on the whole
complex plane, and satisfies the functional equation
$$m(z)+\bar{m}(\bar{z})=-2\sum_{n\ge1}\frac{\mu(n)}{n}\cos(\frac{2\pi}{n}e^{-z})=A(z).$$
(Here the bar designates the complex conjugate.) In the next section we will offer a criteria for the Riemann hypothesis for a
Laplace transform involving
$$A(-\frac{1}{2}\log(z))=-2\sum_{n\ge1}\frac{\mu(n)}{n}\cos(\frac{2\pi}{n}\sqrt{z}).$$
\section{main integrals}
An explicit formula (which is known [11]) for the Riesz function that will become central in this section is given by $x>0$
\begin{equation}\dot{R}(x)=\sum_{\rho}\frac{x^{\rho/2}\Gamma(1-\frac{\rho}{2})}{\zeta'(\rho)}-\sum_{n\ge1}\frac{n!x^{-n}}{\zeta'(-2n)}.\end{equation}
To see this formula, one may apply the Residue theorem (using a positively oriented circle of radius $\frac{1}{2}+M$ centered at the origin) to
$$g(z)=\frac{x^{z/2}\Gamma(1-\frac{z}{2})}{\zeta(z)},$$
after noting simple poles at the complex zeros of the Riemann zeta function, $\rho=\frac{1}{2}+i\gamma,$ $\gamma\in\mathbb{R},$ the trivial zeros at $z=-2j,$ $j\in\mathbb{N},$ and the
zeros of the gamma function at $z=2i+2,$ $i\in\mathbb{N}_0.$ Further, using the formula [1]
$$\frac{1}{\zeta'(-2n)}=\frac{\pi^{2n}2^{2n+1}}{(-1)^n\zeta(2n+1)(2n)!},$$ for positive integers $n,$
we may write the series on the far right hand side of (2.1) as
\begin{equation}\sum_{n\ge1}\frac{n!x^{-n}}{\zeta'(-2n)}=\frac{1}{2}\sum_{n\ge1}\frac{(-1)^nn!}{\zeta(2n+1)(2n)!}\left(\frac{2\pi}{\sqrt{x}}\right)^{2n}.\end{equation}
The series (2.2) might be realized as a Fourier cosine integral in the following way. First, note that
$$\int_{0}^{\infty}xe^{-ax^2}\cos(xb)dx=\frac{1}{2}\int_{0}^{\infty}e^{-ax}\cos(\sqrt{x}b)dx,$$
and that
\begin{equation}\int_{0}^{\infty}e^{-ax}\cos(\sqrt{x}b)dx=\frac{1}{a}\sum_{n\ge0}\frac{n!}{(2n)!}\left(-\frac{b^2}{a}\right)^n.\end{equation}
Consequently, by uniform convergence, we may write the series in (2.2) as
\begin{equation}\int_{0}^{\infty}t\left(x\sum_{n\ge1}n\mu(n)e^{-x(nt)^2}\right)\cos(t2\pi)dt.\end{equation}
We may also use (2.3) to write this as
\begin{equation}\frac{x}{(2\pi)^2}\int_{0}^{\infty}e^{-xt/(2\pi)^2}\left(\sum_{n\ge1}\frac{\mu(n)}{n}\cos(\frac{\sqrt{t}}{n})\right)dt,\end{equation}
which involves the function of Bartz [2, eq.(2.8)].
\begin{theorem} A necessary and sufficient condition for the Riemann hypothesis, is (for all $\epsilon>0$)
\begin{equation}\sum_{\rho}\frac{x^{\rho/2}\Gamma(1-\frac{\rho}{2})}{\zeta'(\rho)}-\frac{x}{(2\pi)^2}\int_{0}^{\infty}e^{-xt/(2\pi)^2}\left(\sum_{n\ge1}\frac{\mu(n)}{n}\cos(\frac{\sqrt{t}}{n})\right)dt=O(x^{\frac{1}{4}+\epsilon}).\end{equation}
\end{theorem}

We now turn our attention to a useful study on Fourier integrals employed by Csordas [4]. The following definition, along with applicable theorems, was used there (see also that papers' references) to determine the nature of the zeros and other properties of the function represented by the Fourier cosine integral
$$f(x):=\int_{0}^{\infty}k(t)\cos(xt)dt.$$
These type of results had its beginnings with that of P\'{o}lya [9]. 
\\*
\\*
{\bf Definition 2.2} \it A function $k:\mathbb{R}\rightarrow\mathbb{R}$ is said to be an \it 'admissible kernel,' if it 
satisfies (i) $k(t)\in C^{\infty}(\mathbb{R}),$ (ii) $k(t)>0$ for $t\in\mathbb{R},$ (iii) $k(t)=k(-t)$ for $t\in\mathbb{R},$ (iv) $\frac{d}{dt}k(t)<0$ for $t>0,$ and (v) for some $\epsilon>0,$ 
$$k^{(n)}(t)=O(e^{-|t|^{2+\epsilon}}),$$
as $t\rightarrow\infty.$ \rm
\\*
Some interesting properties are known about $f(x)$ when $k(t)$ satisfies Definition 2.2 [4]. Namely, by the Riemann-Lebesgue Lemma,
$f(x)\rightarrow0$ as $|x|\rightarrow\infty.$ Additionally, $f(x)$ is then an entire function of order $\frac{2+\epsilon}{1+\epsilon}<2.$
\\*
\par 
We may observe that we may make the the change of variables in (2.4) with $t$ replaced by $t/\sqrt{x}$ 
to get that
\begin{equation}\int_{0}^{\infty}t\left(\sum_{n\ge1}n\mu(n)e^{-(nt)^2}\right)\cos(t2\pi/\sqrt{x})dt,\end{equation}
still represents the series in (2.2). However, our $k(t)$ function here is odd and subsequently it is not the case that the function $f(x)$ represented by our integral can have only real zeros. Our $k(t)$ implies that our $f(x)$ has infinitely many non-real zeros, and finitely many real zeros. To see this, we need only observe that if we let $s(t)=\sum_{n\ge1}n\mu(n)e^{-(nt)^2},$ then $s(t)>0$ when $t>0,$ $s(t)=s(-t),$ and $s'(t)<0$ when $t>0.$ Comparing these properties with the work of Csordas [4] gives our claim. In fact, it was already shown in Rieszs' study [10] that $\dot{R}(x)$ has infinitely many imaginary zeros using a different approach. See also [13] for more properties on the zeros of the Riesz function.
\section{Remarks on a recent generalization}
 Recently, Dixit et. al. [5, 6] studied a more general series than a function considered by Hardy and Littlewood, with a similar condition to (1.1). Their function is given by
\begin{equation} P_{z}(y)=\sum_{n\ge1}\frac{\mu(n)}{n}e^{-y/n^2}\cosh(\frac{\sqrt{y}z}{n}).\end{equation} They
offer the condition that the Riemann hypothesis implies $P_{z}(y)=O_{z,\epsilon}(y^{-\frac{1}{4}+\epsilon})$ as $y\rightarrow\infty$ for
all $\epsilon>0.$ We offer some more comments on $P_{z}(y)$ herein. It is well-known [7] that (for $y>0$ and $\Re(\beta)>0$)
\begin{equation} \int_{0}^{\infty}e^{-t^2/{4\beta}}\cosh(\alpha t)\cos(yt)dt=\sqrt{\pi/\beta}e^{\alpha^2\beta}e^{-\beta y^2}\cos(2\alpha\beta y).\end{equation}
Or equivalently,
\begin{equation} e^{-t^2/{4\beta}}\cosh(\alpha t)=\sqrt{\pi/\beta}e^{\alpha^2\beta}\int_{0}^{\infty}e^{-\beta y^2}\cos(2\alpha\beta y)\cos(yt)dy.\end{equation}
Put $\beta=n^2,$ and $\alpha=1/n,$ where $n\in\mathbb{N}.$ Then sum over $\mu(n)/n$ to get, by uniform convergence,
\begin{equation} \sum_{n\ge1}\frac{\mu(n)}{n}e^{-t^2/{4n^2}}\cosh(\frac{t}{n})=\sqrt{\pi}e\int_{0}^{\infty}\sum_{n\ge1}\frac{\mu(n)}{n^2}e^{-n^2 y^2}\cos(2n y)\cos(yt)dy.\end{equation}
We need to check that the function
$$\bar{k}(t)=\sum_{n\ge1}\frac{\mu(n)}{n^2}e^{-n^2 t^2}\cos(2n t),$$
is an admissible kernel according to Definition 2.2. The condition (ii) that $\bar{k}(t)>0$ for $t\in\mathbb{R},$ can not be met for all $t\in\mathbb{R}.$ 
For example, in considering $\bar{k}(t)$ partial sums, we see that $t=\pi/2$ gives
\begin{equation}\sum_{1\le n \le 3}\frac{\mu(n)e^{-n^2\pi^2/4}\cos(n\pi)}{n^2}=-\frac{1}{36}e^{-\frac{9\pi^2}{4}}\left(-4+9e^{5\pi^2/4}+36e^{2\pi^2}\right), \end{equation}
which is $<0.$ Subsequently (ii) only holds for a subset of $\mathbb{R}.$ If instead we choose $\alpha=i/n,$ initially in our computations from (3.3), we find that applying our same analysis leads to an admissible kernel for the function $P_{iz'}(y),$ $z'\in\mathbb{R},$ which implies it would have only real zeros if $z'\in\mathbb{R}.$
\par Riesz [10] noted that $|\dot{R}(x)|<|x|e^{|x|},$ since $x\in\mathbb{R},$ and hence $\dot{R}(x)$ has order one. \par Recall the Hermite polynomials are generated by
$$e^{2xt-t^2}=\sum_{n\ge0}\frac{H_n(x)t^n}{n!},$$ and the Hermite numbers are $H_n:=H_n(0).$
\begin{theorem} The function $\dot{R}(x^2)/x^2$ may take the form
$$\sum_{n\ge0}\frac{H_nx^n}{n!\zeta(n+2)},$$
where $H_n$ is the $n$th Hermite number, and therefore is an entire function of order $\lambda,$ given by
$$\lambda=\lim_{n\rightarrow\infty}\sup\frac{n\log{n}}{\log{\frac{|n!\zeta(n+2)|}{|H_n|}}}.$$\end{theorem}
\section{More on the Riesz Criterion and other possible directions}
We make some further comments on Theorem 2.1 and offer some possible further directions. First if we first note that
the series on the far right hand side of (2.1) is $O(\frac{1}{x}),$ we may write
\begin{equation}\dot{R}(x)=\sum_{\rho}\frac{x^{\rho/2}\Gamma(1-\frac{\rho}{2})}{\zeta'(\rho)}+O(\frac{1}{x}).\end{equation}
Let
\begin{equation}\sum_{\rho}\frac{x^{\rho/2}\Gamma(1-\frac{\rho}{2})}{\zeta'(\rho)}=O(w(x)).\end{equation}
Then
\begin{equation}\dot{R}(x)=O(\max\{w(x),\frac{1}{x}\}).\end{equation}
This, together with Riesz criterion (1.1), tells us that proving $w(x)=x^{\frac{1}{4}+\epsilon},$ $\epsilon>0,$ would imply the Riemann Hypothesis. This is equivalent to the observation that (4.1) says $\dot{R}(x)\sim \sum_{\rho}\frac{x^{\rho/2}\Gamma(1-\frac{\rho}{2})}{\zeta'(\rho)}.$
\par A possible direction for further research would be to consider the following integral, which we produce
from some observations. We start with Riesz's integral [10]
\begin{equation}\frac{\Gamma(1-\frac{s}{2})}{\zeta(s)}=\int_{0}^{\infty}x^{-(\frac{s}{2}+1)}\dot{R}(x)dx.\end{equation}
Now if we assume a suitable test function $T(x)$ has a Mellin transform $\bar{T}(s),$ which exists in the region $\frac{1}{2}+\eta\le\Re(s)\le 2-\eta,$ ($\eta>0$), then 
\begin{equation}\int_{0}^{\infty}x^{s-1}\left(\int_{0}^{\infty}T(x\sqrt{t})\frac{\dot{R}(t)}{t}dt\right)dx=\frac{\bar{T}(s)\Gamma(1-\frac{s}{2})}{\zeta(s)}.\end{equation}
So we may write 
\begin{equation}\int_{0}^{\infty}T(x\sqrt{t})\frac{\dot{R}(t)}{t}dt=\frac{1}{2\pi i}\int_{c-i\infty}^{c+i\infty}\frac{\bar{T}(s)\Gamma(1-\frac{s}{2})}{\zeta(s)}x^{-s}ds.\end{equation}
If we choose our test function to be the dirac delta function $T(t)=\delta(t-a)$ and choosing $c=\frac{1}{2}+\epsilon,$ we can obtain the Riesz condition after noting that Littlewood [8, pg.161] showed that a criterion for the Riemann Hypothesis is that $\sum_{n\ge1}\mu(n)n^{-\frac{1}{2}-\epsilon}$ converges for every $\epsilon>0.$ It would be interesting to see some further examples by choosing other test functions $T(x).$ For example, choosing a function $T_{\epsilon_1}(t)$ whose limit is $\lim_{\epsilon_1\rightarrow0}T_{\epsilon_1}(t)=T(t)=\delta(t-a),$ and estimating the integral in (4.6) involving $T_{\epsilon_1}(t)$ may lead to improvements or new criteria.

{\bf Acknowledgement.} We thank Professor Dixit and Professor Wolf for helpful comments.

1390 Bumps River Rd. \\*
Centerville, MA
02632 \\*
USA \\*
E-mail: alexpatk@hotmail.com

\begin{thebibliography}{9}
\bibitem{ConcreteMath} 
G. Andrews, R. Askey, and R. Roy. \emph{Special Functions,} volume 71 of Encyclopedia of Mathematics and its Applications. Cambridge University Press, New York, 1999.
\bibitem{ConcreteMath}
K. M. Bartz, \emph{On some complex explicit formulae connected with the M$\ddot{o}$bius function, I,} Acta Arithmetica 57, p.283--293 (1991).
\bibitem{ConcreteMath}
K. M. Bartz, \emph{On some complex explicit formulae connected with the M$\ddot{o}$bius function, II,} Acta Arithmetica 57, p.295--305 (1991).
\bibitem{ConcreteMath}
G. Csordas, \emph{Fourier Transforms of Positive Definite Kernels and the Riemann $\xi$-function,} Computational Methods and Function Theory, Volume 15, Issue 3, pp 373--391 (2015).
\bibitem{ConcreteMath}
A. Dixit, \emph{Analogues of the general theta transformation formula,} Proc. Roy. Soc. Edinburgh, Sect. A, 143 (2013), 371--399
\bibitem{ConcreteMath}
A. Dixit, A. Roy and A. Zaharescu, \emph{Riesz-type criteria and theta transformation analogues,} J. Number Theory 160, p. 385--408 (2016).
\bibitem{ConcreteMath} 
I. S. Gradshteyn and I. M. Ryzhik, eds., Table of Integrals, Series, and Products, 7th ed., Academic Press,
San Diego, 2007.
\bibitem{ConcreteMath} 
G. H. Hardy and J. E. Littlewood, \emph{Contributions to the Theory of the Riemann Zeta-Function and the
Theory of the Distribution of Primes,} Acta Math., 41 (1916), 119--196.
\bibitem{ConcreteMath}
G. P\'{o}lya, \emph{Uber trigonometrische Integrale mit nur reellen Nullstellen,} J. Reine Angew. Math. 158 (1927), 6--18.
\bibitem{ConcreteMath} 
M. Riesz, \emph{Sur l'hypoth'ese de Riemann,} Acta Math., 40 (1916), 185--190.
\bibitem{ConcreteMath}
G.W. Smith, \emph{On a function of Marcel Riesz,} http://arxiv.org/abs/1209.5652, September 2012.
\bibitem{ConcreteMath} 
E. C. Titchmarsh, \emph{The Theory of the Riemann Zeta Function,} Clarendon Press, Oxford, 1986.
\bibitem{ConcreteMath}
H.Wilf, \emph{On the zeros of Riesz' function in the analytic theory of numbers,} Illinois J. Math., 8 (1964), pp. 639--641
\bibitem{ConcreteMath}
M. Wolf, \emph{Evidence in favor of the Baez-Duarte criterion for the Riemann Hypothesis,} 
Computational Methods in Science and Technology, v.14 (2008) pp.47--54

\end{thebibliography}
\end{document}